\documentclass[12pt]{amsart}
\usepackage{latexsym,amssymb,amsmath}
\usepackage[all,ps,cmtip]{xy} 
\newtheorem{theo}{Theorem}
\newtheorem{coro}{Corollary}

\newtheorem{lemm}{Lemma}

\begin{document}

\def\ot{\otimes}
\def\we{\wedge}
\def\wec{\wedge\cdots\wedge}
\def\op{\oplus}
\def\ra{\rightarrow}
\def\lra{\longrightarrow}
\def\fso{\mathfrak so}
\def\cO{\mathcal{O}}
\def\cS{\mathcal{S}}
\def\fsl{\mathfrak sl}
\def\PP{\mathbb P}\def\PP{\mathbb P}\def\ZZ{\mathbb Z}\def\CC{\mathbb C}
\def\RR{\mathbb R}\def\HH{\mathbb H}\def\OO{\mathbb O}
\def\smc{\cdots }
\title{A note on certain Kronecker coefficients}

\author{L. Manivel}

\begin{abstract} 
We prove an explicit formula for the tensor product with itself of an irreducible complex 
representation of the symmetric group defined by a rectangle of height two. We also describe
part of the decomposition for the tensor product of representations defined by rectangles 
of heights two and four. Our results are deduced, through Schur-Weyl duality, from the observation
that certain actions on triple tensor products of vector spaces, are multiplicity free. 
\end{abstract}

\maketitle


\section{Introduction}
Irreducible complex representations of the symmetric group $\cS_n$ are well known to be indexed 
by partitions of $n$ in a natural way (see e.g. \cite{mcd}, I.7). We will denote by $[\lambda]$ the representation 
associated to the partition $\lambda$. A major unsolved problem is to find a general rule for the tensor 
product $[\lambda]\otimes [\mu]$ of two such representations. Equivalently, one would like a general 
rule for the computation of the Kronecker coefficients, which are defined as the multiplicities appearing 
in the formula
$$[\lambda]\otimes [\mu] = \bigoplus_{\nu} k_{\lambda\mu\nu}[\nu].$$

In the setting of algebraic complexity theory, a specific instance of this general problem has been put to
the fore: can one compute the tensor product with itself, of an irreducible complex 
representation of the symmetric group defined by a {\it rectangle} partition? (see \cite{bklmw}
for an overview). If the rectangle 
has height one this is pretty obvious, since the corresponding representation is the trivial one. 
In this note we give an answer for the next case, that of a rectangle of height two. (Note that 
a rectangle of width two would lead exactly to the same answer, since one can pass from a partition 
to the dual one, in terms of representations of the 
symmetric group, simply by the product with the sign representation.)

Our main result is the following, where the length $\ell(\lambda)$ of a partition $\lambda$ is 
the number of its non zero parts :

\begin{theo}
For any integer $n$, the tensor product $[n,n]\otimes [n,n]$ is multiplicity free.
Its decomposition is 
$$[n,n]\otimes [n,n] = \bigoplus_{\substack{\lambda\,\mathrm{even},\; |\lambda|=2n \\ \ell(\lambda)\le 4}}[\lambda]
\oplus\bigoplus_{\substack{\mu\,\mathrm{odd},\; |\mu|=2n\\ \ell(\mu)=4}}[\mu].$$
\end{theo}

It would be interesting to understand the splitting of $[n,n]\otimes [n,n]$
into its symmetric and skew-symmetric parts. 

\medskip
To state our second result, we introduce the following notation: 
$$[\lambda]\otimes_\ell [\mu] = \bigoplus_{\ell(\nu)\le\ell} k_{\lambda\mu\nu}[\nu].$$

\begin{theo}
For any integer $n$, the partial tensor product $[2n,2n]\otimes_3 [n,n,n,n]$ is multiplicity free.
Its decomposition is 
$$[2n,2n]\otimes_3 [n,n,n,n] = \bigoplus_{\substack{|\lambda|=2n \\ 
\lambda_2+\lambda_3-\lambda_1\ge 0\,\mathrm{and\,even}}}[2\lambda].$$
\end{theo}

\section{Schur-Weyl duality and multiplicity free actions}
In order to prove the previous two theorems we will restate them in terms of 
representations of general linear groups, in a quite standard way. Recall the statement 
of the Schur-Weyl duality between representations of symmetric groups and of general 
linear groups: let $V$ be any finite dimensional complex vector space, and $n$ any integer. 
Then the $\cS_n\times GL(V)$-module $V^{\otimes n}$ decomposes as 
$$V^{\otimes n} = \bigoplus_{|\lambda|=n}[\lambda]\otimes S_\lambda V,$$
where $S_\lambda V$ denotes the Schur module of weight $\lambda$, which is an irreducible 
polynomial representation of $GL(V)$.
A straightforward consequence is that, for three vector spaces $U,V,W$ and three partitions 
$\lambda,\mu,\nu$ of the same integer $n$, the multiplicity of $S_\lambda U\otimes S_\mu V
\otimes S_\nu W$ inside $S^n(U\otimes V\otimes W)$ is equal to the Kronecker coefficient 
$k_{\lambda\mu\nu}$. In particular, if $U$ and $V$ have respective dimensions $du$ and  $dv$, 
with $(u,v)=1$, we deduce that
$$Sym(U\otimes V\otimes W)^{SL(U)\times SL(V)}=\bigoplus_{n\ge 0}\bigoplus_{|\lambda|=nuv}
k_{(nv)^u,(nu)^v,\lambda}S_\lambda W.$$
Here $(nv)^u$ denotes the rectangular partition with $u$ parts all equal to $nv$, so that 
the corresponding Schur module of $U$ consists of $SL(U)$-invariants. This shows that 
Kronecker coefficients involving rectangular partitions are closely related to invariant 
theory. Indeed our two theorems above will be translated into the statements that two 
invariant algebras $Sym(U\otimes V\otimes W)^{SL(U)\times SL(V)\times N}$ are polynomial
algebras, where $N$ is a group of strictly upper triangular matrices in $SL(W)$. 

Note that the complete invariant algebra $Sym(U\otimes V\otimes W)^{SL(U)\times SL(V)\times SL(W)}$
is then also a polynomial algebra. This happens when the dimensions of the three
spaces are either $(n,2,2)$, $(2,3,3)$,  $(2,3,4)$ or  $(2,3,5)$. In the terminology of \cite{kac},
these cases correspond to $\theta$-groups defined by the triple nodes of the Dynkin diagrams
of type $D_{n+2}$, $E_6$, $E_7$ and $E_8$. The two cases we examine in this note are thus related to
$D_6$ and $E_7$, respectively. 

\smallskip
As is well known, multiplicity free actions of reductive groups can be detected by the existence 
of an open orbit for a Borel subgroup. We use this principle in the following setting: let 
$G$ and $H$ be two reductive groups with finite dimensional representations $V$ and $W$. 
Let $B$ denote a Borel subgroup of $H$ and $N$ its unipotent radical.  
Suppose that $G\times B$ acts on $V\otimes W$
with an open orbit $\cO$. Let $X_1,\ldots ,X_r$ denote the boundary components of $\cO$, that
is, the irreducible hypersurfaces in its complement. Applying \cite{brion}, Proposition 3 of Chapter 3, 
we are lead to the following conclusions:
\begin{itemize}
\item $X_1,\ldots ,X_r$ have equations $f_1,\ldots ,f_r$ which are semi-invariants of $B$ 
with linearly independant weights $\mu_1,\ldots ,\mu_r$; in particular $r$ cannot exceed the 
rank of $H$;  
\item the algebra $\CC[V\otimes W]^{G\times N}$ is a polynomial algebra over $f_1,\ldots ,f_r$. 
\end{itemize}
That $f_1,\ldots ,f_r$ are semi-invariants of $B$ of weights $\mu_1,\ldots ,\mu_r$ means that 
$f_i(bx)=\mu_i(b)f_i(x)$ for all $x\in V\otimes W$ and $b\in B$. 
Moreover, 
as an $H$-module, the algebra of $G$-invariant functions on $V\otimes W$ is multiplicity-free:
$$\CC[V\otimes W]^G = \bigoplus_{\mu\in\ZZ_+\mu_1+\cdots +\ZZ_+\mu_r}W_\mu,$$ 
if $W_\mu$ denotes the irreducible $H$-module of highest weight $\mu$. Indeed, such a component of 
$\CC[V\otimes W]^G$ can be detected by its one-dimensional subspace of $N$-invariants. 

\smallskip
In order to prove our two theorems, we will therefore just need to prove that the corresponding actions have 
open orbits, and to identify the boundary components.

\section{Proof of Theorem 1}

Let $U$ and $V$ be two-dimensional vector spaces.  

\begin{lemm} 
Consider the action of $SL(U)\times SL(V)$ on the flag variety $\mathcal{F}(U\otimes V)$.
The generic isotropy group of this action is a product of $\ZZ_2$ by a quaternion group.   
\end{lemm}
 
\proof
Consider a general flag $W_1\subset W_2\subset W_3\subset U\otimes V$. The projective 
line $\PP W_2\subset \PP (U\times V)$ meets the quadric $Q=\PP U\times \PP V$ in two general
points, which means that $W_2$ has a basis  of the form $u_0\otimes v_0, u_1\otimes v_1$, where 
$u_0,u_1$ is a basis of $U$ and $v_0,v_1$ is a basis of $V$.  Multiplying if necessary, one of these
vectors by a scalar, we may suppose that $W_1$ is the line in $W_2$ generated by $u_0\otimes v_0+
u_1\otimes v_1$. Finally, $W_3$ is the kernel of a general linear form $\phi$ vanishing on $W_2$. 
Since in terms of the dual basis, $W_2^\perp=\langle u_0^\vee\otimes v_1^\vee, u_1^\vee\otimes v_0^\vee
\rangle$, 
we can suppose that $\phi=u_0^\vee\otimes v_1^\vee - u_1^\vee\otimes v_0^\vee$. This means that $W_3$ is 
generated by $W_2$ and $u_0\otimes v_1+u_1\otimes v_0$. 

Now it is straightforward to compute the stabilizer of our flag explicitely, and to identify it
with the product of $\ZZ_2$ by a quaternion group. \qed

\medskip
In fact the only important thing to us is that this stabilizer is finite, because of the
following corollary.  
Let $W$ be a four-dimensional vector space, and $B$ a Borel subgroup in $GL(W)$. 

\begin{coro}
The group $SL(U)\times SL(V)\times B$ has an open orbit in $U\otimes V\otimes W$. 
\end{coro}

\proof 
Consider a tensor $T\in U\otimes V\otimes W$ as a morphism $\phi_T : W^\vee\rightarrow U\otimes V$. 
For a generic $T$ this morphism is injective and maps the flag defining $B$ (or rather the orthogonal
flag) to a generic flag in  $U\otimes V$. By Lemma 1, $SL(U)\times SL(V)$ has an open orbit in the 
flag variety $\mathcal{F}(U\otimes V)$. And once the image flag is fixed, it is clear that $B$ acts 
transitively on the set of compatible injections.  
\qed

\medskip
As we explained above, the next step is to describe the boundary components of the open orbit. 
Let us denote by $F=(W_1\subset W_2\subset W_3\subset W)$ the flag whose stabilizer is 
the Borel subgroup $B$ of $GL(W)$, and the orthogonal flag  in  $W^\vee$ by  $F^\perp$. 

As in the proof of the Lemma we denote by $\phi_T : W^\vee\rightarrow U\otimes V$ the morphism 
defined by the tensor $T\in U\otimes V\otimes W$. We can describe the boundary components of the 
open orbit in $U\otimes V\otimes W$ by the following codimension one conditions:
\begin{enumerate}
\item 
$\phi_T$ is not an isomorphism. The corresponding boundary component $X_1$ is the complement of the 
$SL(U)\times SL(V)\times GL(W)$-orbit. It is just the quartic hypersurface of equation $f_1 =\det 
\phi_T$. 
This equation is a weight vector in $S_{22}U^\vee\otimes S_{22}V^\vee\otimes \wedge^4W^\vee=
\wedge^4(U\otimes V)^\vee\otimes\wedge^4W^\vee\subset S^4(U\otimes V\otimes W)^\vee$. This means that
the weight $\mu_1$ of $f_1$ is, written as a sequence of three partitions, $\mu_1=(22,22,1111)$. 
\item 
$\phi_T(\PP W_3^\perp)$ belongs to the quadric $Q$. The corresponding boundary component $X_2$
is defined by the condition that $q(\phi_T(w^\vee))=0$, if $w^\vee$ generates $W_3^\perp$ and $q$ denotes an
equation of $Q$. Thus an equation $f_2$ of $X_2$ is a highest weight vector in 
$\wedge^2U^\vee\otimes \wedge^2V^\vee\otimes S_2W^\vee$.
It has degree two and weight $\mu_2=(11,11,2)$. 
\item 
$\phi_T(\PP W_2^\perp)$ is a tangent line to $Q$.
This is the case if and only if $\phi_T(W_2^\perp)$ is generated
by vectors of the form $u_0\otimes v_0$ and $u_0\otimes v_1+u_1\otimes v_0$. Considered as a line $\ell$ in 
$\wedge^2W^\vee$, this means that $W_2^\perp$  is mapped by $\phi_T$ to the line generated by  
$u_0^2\otimes (v_0\wedge v_1)\oplus (u_0\wedge u_1)\otimes v_0^2$ in $\wedge^2(U\otimes V)=
S^2U\otimes \wedge^2V\oplus \wedge^2U\otimes S^2V$. 
Since the map $S^2(S^2U)\rightarrow S_{22}U$ kills any tensor of the form $(u^2)^2$, we deduce that $\phi_T$
maps $\ell^2$ to zero in $S_{22}U\otimes S_{22}V$. This implies that an equation $f_3$ of the corresponding
boundary component $X_3$ is a highest weight vector in $S_{22}U^\vee\otimes S_{22}V^\vee\otimes S_{22}W^\vee$.
It has degree four and its weight is $\mu_3=(22,22,22)$. 
\item
$\phi(\PP W_1^\perp)$ is a tangent plane to $Q$. 
This is similar to the case of $X_2$, up to duality. A hyperplane $H$ in $U\otimes V$ defines a line in 
$\wedge^3(U\otimes V)=U\otimes V\otimes (\wedge^2U\otimes \wedge^2V)$, hence a line $\ell$ in $U\otimes V$
(the orthogonal line with respect to the polarity defined by $Q$). 
This hyperplane is tangent to the quadric $Q$ if and only if $\ell$ is contained in $Q$. This means that
$X_4$ is defined by the condition that the composition 
$$S^2(\wedge^3W^\vee)\rightarrow S^2(\wedge^3(U\otimes V))=S^2(U\otimes V)\otimes (\wedge^2U\otimes \wedge^2V)^2
\rightarrow (\wedge^2U\otimes \wedge^2V)^3$$
vanishes. Hence  an equation $f_4$ of $X_4$ is a highest weight vector in 
$S_{33}U^\vee\otimes S_{33}V^\vee\otimes S_{222}W^\vee$. It has degree six and weight 
$\mu_4=(33,33,222)$.
\end{enumerate}
The four weights of $f_1,f_2,f_3,f_4$ are linearly independent. Since the rank of $GL(W)$ is four, 
we must have found all the boundary components and we can conclude that 
$$\CC[U\otimes V\otimes W]^{SL(U)\times SL(V)\times N}=\CC[f_1,f_2,f_3,f_4],$$
where $N$ denotes the unipotent radical of $B$. 
This implies that $\CC[U\otimes V\otimes W]$ contains a copy of 
$S_{n,n}U^\vee\otimes S_{n,n}V^\vee\otimes S_{\lambda}W^\vee$ if and only if $\lambda$ is a non negative
linear combination of the components of $\mu_1,\mu_2,\mu_3$ and $\mu_4$ on $W$, that is, the weights 
$(1111), (2), (22)$ and $(222)$. Moreover, in that case the multiplicity 
is equal to one. 

Rephrasing this result via Schur-Weyl duality we get the statement of Theorem 1. 

\section{Proof of Theorem 2}


For the proof of the next Lemma we need to recall  briefly
the principle of {\it castling transforms}, introduced by Sato and Kimura \cite{sk}.
Consider a $G$-module $M$
of dimension $m$, and a vector space $N$ of dimension $n$. Suppose that $m>n$. A tensor $T$ in $M\otimes N$ 
can be identified to a linear map $\phi_T: N^\vee\rightarrow M$. If $\phi_T$ is injective, in particular
for a generic $T$, the stabilizer of $T$ in $G\times GL(N)$ is canonically isomorphic to the stabilizer 
in $G$ of the image of $\phi_T$, considered as a point of the Grassmannian $G(n,M)$. But this Grassmannian 
is isomorphic with $G(n-m,M^\vee)$, and the generic stabilizer of the action of $G\times GL(N)$ 
on $M\otimes N$ is therefore isomorphic with the generic stabilizer of the action of $G\times GL(P)$
on $M^\vee\otimes P$, for $P$ a vector space of dimension $m-n$. Replacing  $M\otimes N$ by $M^\vee\otimes P$
is precisely what Sato and Kimura call a castling transform. In case $M\otimes N$ is prehomogeneous
and $m<2n$, $M^\vee\otimes P$ is also prehomogeneous but of smaller dimension. 

\medskip
Let $U,V,W$ be complex vector spaces of respective dimension two, 
four and three. Let $B$ be a Borel subgroup of $GL(W)$. 

\begin{lemm}
The group $SL(U)\times SL(V)\times B$ has an open orbit in 
$U\otimes V\otimes W$. 
\end{lemm}

\proof We claim that the generic isotropy group of the action of $SL(U)\times SL(V)\times GL(W)$ 
on $U\otimes V\otimes W$ is a copy of $SL_2$, up to a finite group. To check this, we oberve that 
$U\otimes V\otimes W$ is, according to the terminology of Sato and Kimura, a non-reduced prehomogeneous 
vector space, which means that it is related to smaller prehomogeneous spaces of the same type by 
certain castling transforms. 

In order to apply this process to the case we are interested in, we first observe that the generic stabilizers
of $SL(U)\times SL(V)\times GL(W)$ and $SL(U)\times GL(V)\times SL(W)$on $U\otimes V\otimes W$ are equal,
up to a finite group. Applying a castling transform with $G=SL(U)\times SL(W)$ acting on $U\otimes W$, we 
deduce that this generic stabilizer is the same as the generic stabilizer 
of  the action of $SL(U)\times GL(Q)\times SL(W)$ on $U^\vee\otimes Q\otimes W^\vee$, where $Q$ has dimension 
$2\times 3-4=2$. 
Up to a finite group, this is also the generic stabilizer 
of  the action of $SL(U)\times SL(Q)\times GL(W)$, and after a new castling trasform, we deduce that 
this is also the generic stabilizer of  the action of $SL(U)\times SL(Q)\times GL(R)$ on 
$U\otimes Q^\vee\otimes R$, where now $R$ has dimension $2\times 2-3=1$. Let us identify $U$ and $Q$, 
which are both two-dimensional. Then it is easy to see that the identity map $I\in U\otimes Q^\vee$ has 
generic stabilizer, so that this generic stabilizer is just 
a copy of $SL_2$ embedded diagonally in $SL(U)\times SL(Q)$.

We can keep track of this generic stabilizer along our two castling transforms. We 
start from the point in $U\otimes U^\vee =End(U)$ defined by $I$. 
The corresponding point in $U^\vee\otimes U\otimes W^\vee$, where $W$ is identified with the orthogonal
to $I$ in $End(U)$ (the hyperplane $End_0(U)$ of traceless matrices), is just the graph of the embedding
of  $End_0(U)$ in $End(U)$. Its isotropy is the image of $SL(U)$ in $SL(U^\vee)\otimes SL(U)\otimes SL(W^\vee)$
given by the natural action of $SL(U)$ on each of the three spaces. Now we make our second castling trasform to get a
point in $U\otimes V\otimes W$, where $V$ is now identified with the kernel of the natural 
evaluation map $End_0(U)\otimes U\rightarrow U$. The corresponding stabilizer is again a copy of 
$SL(U)$ embedded in $SL(U)\otimes SL(V)\otimes SL(W)$ through its natural action on $U,V$ and $W$. 

We can now check our claim: it can be translated into the assertion that the image of $SL(U)$ into
$GL(W)$, where $W=End_0(U)$, does not intersect a general Borel subgroup. But this is straightforward:
such a Borel subgroup is defined by a line, generated by a generic traceless matrix $m$, and a hyperplane
containing it, which can be defined as the orthogonal to a generic traceless matrix $n$ orthogonal to $m$. 
For an element of $SL_2$, preserving $m$  and $n$ forces it to belong to the intersection of two 
tori, and this intersection is finite. \qed

\medskip
Now we identify in $U\otimes V\otimes W$ 
the boundary components of the open orbit of $SL(U)\times SL(V)\times B$.
They can be described in terms of the following codimension one conditions on a tensor 
$T\in U\otimes V\otimes W$, which will be best expressed in terms of certain auxiliary morphisms. 
\begin{enumerate}
\item 
First recall that $SL(U)\times SL(V)\times GL(W)$ has itself an open orbit whose complement is an irreducible 
hypersurface $X_1$ of degree $12$ \cite{kac}. An equation $f_1$ of this hypersurface can be obtained as follows. 
Consider the morphism $\psi^T : U^\vee\otimes W^\vee\rightarrow V$ induced by $T$. Taking is second wedge power, 
we get an induced map 
$$\Psi^T : \wedge^2U^\vee \otimes S^2W^\vee\hookrightarrow 
\wedge^2(U^\vee \otimes W^\vee)\stackrel{\wedge^2\psi^V}{\longrightarrow}\wedge^2V$$
between two vector spaces of the same dimension, six. We can thus let $f_1=\det\Psi^T$, an
invariant of degree $12$ and weight $\mu_1=(66,3333,444)$. 
\item 
Restricting $\psi^T$, we can define a morphism 
$$\psi^T_1 : U\otimes W_1^\perp\rightarrow V$$ 
between two vector spaces of the same dimension, four. We can thus define another boundary component $X_2$
by the condition that this is not an isomorphism. An equation  of this hypersurface is 
$f_2=\det\psi^T_1$, a semi-invariant of degree $4$ and weight $\mu_2=(22,1111,22)$. 
\item
To describe our next boundary component, we need to recall that there exists an invariant non-degenerate
skew-symmetric form of $S^3U$, or equivalently an equivariant morphism 
$$\omega:\wedge^2(S^3U)\rightarrow (\wedge^2 U)^3.$$ 
Now consider the morphism
$$\phi_T :\wedge^3W^\vee\longrightarrow S^3U\otimes  \wedge^3V\simeq S^3U\otimes  \wedge^4V\otimes V^\vee$$
induced by $T$. (Note that $\wedge^3W$ and $\wedge^4V$ are both one dimensional.) 
Taking its square, we get a map 
$$\Phi_T : S^2(\wedge^3W^\vee)\rightarrow \wedge^2(S^3U)\otimes (\wedge^4V)^2\otimes \wedge^2V^\vee
\rightarrow (\wedge^2 U)^3\otimes (\wedge^4 V)^2\otimes \wedge^2V^\vee,$$
where the last arrow is induced by $\omega$. The image of $\Phi_T$ defines, up to scalar, a skew-symmetric form 
$\omega_T$ on $V$. On the other hand, we can restrict the morphism $\Psi^T$ to the line 
$\wedge^2U^\vee \otimes S^2W_2^\perp \subset \wedge^2U^\vee \otimes S^2W^\vee$. Its image is, up to scalar, 
an element of $\Omega_T$ of $\wedge^2V$. We can thus define a boundary component $X_3$ by the condition 
that the natural pairing $\langle\omega_T,\Omega_T\rangle =0$. An equation $f_3$ of this hypersurface has degree 
$8$ and weight $\mu_3=(44,2222,422)$.  
\end{enumerate}
Since we have find three boundary components and $\dim (U\otimes V\otimes W)-\dim (SL(U)\times SL(V)\times B)=3$, 
we must have found all the boundary components and we can conclude that 
$$\CC[U\otimes V\otimes W]^{SL(U)\times SL(V)\times B}=\CC[f_1,f_2,f_3].$$
The weights of $f_1,f_2,f_3$ are independent, as expected, and we deduce that $\CC[U\otimes V\otimes W]$
contains a component $S_{n,n}U\otimes S_{n,n}V\otimes S_{\lambda}W$ if and only if 
$\lambda=(\lambda_1,\lambda_2,\lambda_3)$ is a non negative
linear combination of the components of $\mu_1,\mu_2$ and $\mu_3$ on $W$, that is, 
$(422),(444)$ and $(422)$. Because of the identity
$$\lambda = \frac{\lambda_1-\lambda_2}{2}(422) + \frac{\lambda_3+\lambda_2-\lambda_1}{4}(444) + 
\frac{\lambda_2-\lambda_3}{2}(22),$$
this is equivalent to the conditions that $\lambda$ be even and that $\lambda_3-\lambda_1-\lambda_2$ be a 
non negative multiple of four. Moreover, in that case the multiplicity is equal to one. 
Rephrasing this result via Schur-Weyl duality we get the statement of Theorem 2.

\end{document}